\RequirePackage{fix-cm}
\documentclass[smallextended]{svjour3}       
\smartqed  
\usepackage{amsmath}
\usepackage{graphicx}
\usepackage{amsfonts}
\usepackage{bigfoot}
\usepackage[font=small,labelfont=bf]{caption} 
\usepackage{booktabs} 
\usepackage{amssymb}
\usepackage{subfigure}
\usepackage{color}
\usepackage[ruled,vlined]{algorithm2e}
\DeclareMathOperator{\dist}{dist}

\DeclareMathOperator{\dom}{dom}
\DeclareMathOperator{\interior}{int}

\begin{document}

\title{Global Complexity Analysis of Inexact Successive Quadratic Approximation methods for Regularized Optimization under Mild Assumptions
}


\author{Wei Peng$^1$         \and
	Hui Zhang$^1$ \and Xiaoya Zhang$^1$ 
}


\institute{Wei Peng
	\\\email{weipeng0098@126.com}           
	\\\\
	Hui Zhang(Corresponding author)
	\\\email{h.zhang1984@163.com}
	\\\\
	Xiaoya Zhang
	\\\email{zhangxiaoya09@nudt.edu.cn}
	\\\\
	$^1$ Department of Mathematics, National University of Defense Technology \\
}

\date{Received: date / Accepted: date}

\titlerunning{Inexact SQA without L-smoothness}
\authorrunning{Wei Peng et al.}
\maketitle

\begin{abstract}
Successive quadratic approximations (SQA) are numerically efficient for minimizing the sum of a smooth function and a convex function. The iteration complexity of inexact SQA methods has been analyzed recently. In this paper, we present an algorithmic framework of inexact SQA methods with four types of line searches, and analyze its global complexity under milder assumptions. First, we show its well-definedness and some decreasing properties. {\color{black}Second, under the quadratic growth condition and a uniform positive lower bound condition on stepsizes, we show that the function value sequence and the iterate sequence are linearly convergent. Moreover, we obtain a $o(1/k)$ complexity without the quadratic growth condition, improving existing $\mathcal{O}(1/k)$ complexity results. At last, we show that a local gradient-Lipschitz-continuity condition could guarantee a uniform positive lower bound for the stepsizes.}

\end{abstract}

\section{Introduction}
A fundamental optimization model is ubiquitous in many fields such as machine learning, signal and image processing, and compressed sensing. Typically, the model is to minimize the sum of a smooth function $f$ and a convex regularizer $g$:
\begin{align}\label{model}
\min_{x\in\mathcal{H}} F(x):=f(x)+g(x),
\end{align}
where $\mathcal{H}$ is a real Hilbert space. In recent years, there has been a great deal of interest in developing algorithms to solve (\ref{model}). A classic framework to solve it is the forward-backward splitting(FBS) method \cite{combettes2005signal}, which could be formulated as follows:
\begin{align}\label{approx}
x_{k+1}=\mathop{\arg\min}_{x\in\mathbb{R}^n} \underbrace{f(x_k)+\langle\nabla f(x_k),x-x_k\rangle+\frac{1}{2\tau_k}\|x-x_k\|^2_2}_{\color{black}p_k(x)}+g(x).
\end{align}
Note that $p_{k}$ could be regarded as a second-order approximation to $f$ around $x_k$. Therefore, it is reasonable to further exploit the second-order information of $f$ like
\begin{align*}
q_k(x):=f(x_k)+\langle\nabla f(x_k),x-x_k\rangle+\frac{1}{2}(x-x_k)^T\nabla^2f(x_k)(x-x_k).
\end{align*}
Actually, this kind of approximation was used to design the proximal Newton method \cite{qi1997preconditioning,wei1996convergence}. If $\nabla^2 f$ enjoys some special structure, the proximal Newton method can be very efficient, as shown in \cite{Hsieh_Dhillon_2011} for solving the $\ell_1$-regularized inverse covariance matrix estimation problem. However, in general cases and for large-scale problems, the storage and computation of $\nabla^2 f(x_k)$ could be prohibitive. Thus, one might pursue some approximations to  $\nabla^2 f$. In practice, we usually seek positive definite matrices $H_k$ to approximate $\nabla^2 f(x_k)$. This motivates the proximal successive quadratic approximation (SQA) method:
\begin{align}\label{hkprob}
x_{k+1}=\mathop{\arg\min}_{x\in\mathbb{R}^n} \langle\nabla f(x_k),x-x_k\rangle+\frac{1}{2}(x-x_k)^TH_k(x-x_k)+g(x).
\end{align}
The scheme above can be viewed as a generalization of the FBS and proximal Newton methods, as it reduces to them by setting $H_k=\frac{1}{\tau_k}I_n$ and  $H_k=\nabla^2f(x_k)$, respectively. In this sense, SQA is also known as the proximal quasi-Newton method or the variable metric forward-backward splitting method.

Numerically, one vital issue is how to solve the subproblems (\ref{hkprob}), whose closed-form solutions are usually hard to obtain. Therefore,  iterative algorithms are needed to find inexact solutions to the subproblems. On this road, the authors of \cite{Byrd_Nocedal_Oztoprak_2016} proposed a proximal gradient-based inexactness condition for solving subproblems inexactly, along with a global convergence result. Later, the authors of \cite{yue2016family} showed superlinear and Q-quadratic convergences (with respect to outer iterations) for a family of inexact SQA methods by a delicate parameter choosing strategy for the proximal gradient-based condition. It should be noted that they employed the Luo-Tseng error bound to replace the strong convexity near optimal points.

However, we observe that the inexactness condition might increase the inner iteration complexity as the outer iteration goes on, which is scarcely taken into account in previous works. This is the first motivation of our study.

Besides the gradient-based inexactness condition, which might make the number of inner iterations hard to estimate, some inexactness conditions based on the decrease of function values appeared. For example,  the authors of \cite{Scheinberg_Tang_2016} proposed an inexact SQA method such that the function values of $Q_k(\cdot)$ (which are modifications of the objective functions in (\ref{hkprob})) decrease  to a given absolute error.  With a proper parameter selection strategy and for arbitrary $\epsilon>0$, they showed that $\mathcal{O}(\frac{1}{\varepsilon}\log{\frac{1}{\varepsilon}})$  inner iterations  is enough to yield  an approximate solution $\bar x^\ast$ satisfying
\begin{align}
F(\bar x^\ast)-\inf F\leq \varepsilon.
\end{align}
To put a uniform upper bound on the number of inner iterations, the authors of \cite{lee2018inexact} presented another inexactness condition, which requires that the function values  of $Q_k(\cdot)$ decrease to some relative error. They showed linear convergence results under the optimal set strongly convexity (OSSC), which is  weaker than the strong convexity assumption.

A common assumption made in these existing inexact SQA methods is that the gradient of the smooth part $f$ is Lipschitz continuous, which may fail in many problems. Therefore, it is natural to ask how inexact SQA methods behave without this assumption. This is the second motivation of our study. Before us, we note that a couple of exact gradient-type methods without this assumption have been investigated recently \cite{Cruz_On_2016,Cruz_On_2016_,bello2018q}. Among them, the author of \cite{salzo2017variable} studied a class of exact SQA methods. The main tool employed in \cite{salzo2017variable} is  the quasi-Fej\'er monotone property and thus put a strong restriction on $\{H_k\}$. For inexact SQA methods, we find the line of thought in  \cite{salzo2017variable} is infeasible because the inexact solving of subproblems brings essential difficulties in analyzing iteration points. Besides, we want to drop the additional restriction on  $\{H_k\}$. To this end, we develop new proof methods to study inexact SQA methods under weaker conditions than the gradient-Lipschitz-continuity assumption, and make the following contributions:

\begin{itemize}
	\item[1.] We present an algorithmic framework of inexact SQA methods with four types of backtracking line search strategies, and show its well-definedness without assuming the the gradient-Lipschitz-continuity property.
	\item[2.] We relax OSSC\cite{lee2018inexact} to the weaker quadratic growth condition to obtain Q-linear {\color{black} convergence of the  function value sequence and  R-linear convergence of the iterate sequence.}
	\item[3.] Without the quadratic growth property, we derive a $o(k^{-1})$ convergence
	of the function value sequence by non-trivially modifying \cite[Lemma 6]{lee2018inexact}, which improves the $\mathcal{O}(k^{-1})$ convergence in most	existing related results. The author of \cite{salzo2017variable} also presented a $o(k^{-1})$ convergence
	result but for exact SQA method and with stronger restrictions on $\{H_k\}_{k\geq 0}$.
	\item[4.]  {\color{black}Finally, we show that a local gradient-Lipschitz-continuity condition could guarantee a uniform positive lower bound for stepsizes when performing backtracking line searches.}
\end{itemize}


The paper is organized as follows. In Section 2, we introduce notations and assumptions.
In Section 3, we present the algorithmic framework of inexact SQA methods with four types of backtracking line search strategies. In Section 4, we analyze the global complexity of the algorithmic framework under mild assumptions. In Section 5, we give a short summary of this paper, along with some discussion for future work.

\section{Notations \& Assumptions}
{\color{black}
For a nonempty closed set $\mathcal{C}\subset\mathbb{R}^n$, we denote the distance from $x$ to $\mathcal{C}$ by $\dist(x,\mathcal{C}):=\inf_{y\in\mathcal{C}}\|x-y\|$.
	The domain of an extended-value function $h:\mathbb{R}^n\rightarrow[-\infty,+\infty]$ is defined as
	$\dom h:=\{x\in\mathbb{R}^n: h(x)<+\infty\}$. We say that $h$ is proper if $h(x)>-\infty$ for every $x$ and $\dom h\neq \emptyset$. The gradient of a differentiable function $f$ is denoted by $\nabla f$. } We say that $\nabla f$ is $L$-Lipschitz continuous on a convex set $S$ if
\begin{align*}
\|\nabla f(x)-\nabla f(y)\|\leq L\|x-y\|,~~\forall x,y\in S.
\end{align*}
We say a sequence $\{x_k\}$ R-linearly converges to $x^\ast$ if
\begin{align}
\limsup_{k\rightarrow+\infty}\|x_k-x^\ast\|^\frac{1}{k}<1.
\end{align}
We modify the subproblem \eqref{hkprob} into the following form:
\begin{align}\label{sqa}
\min_{x\in\mathcal{H}} Q_k(x):=\langle \nabla f(x_k),x-x_k\rangle+g(x)-g(x_k)+\frac{1}{2}\|x-x_k\|_k^2,
\end{align}
and let $Q_k^\ast$ be the minimum. The notation $\|\cdot\|_k$ will be explained below.  We denote by $\bar x_{k+1}$ the $\eta$-approximate minimizer, which satisfies
\begin{align}\label{inexact}
Q_k(\bar x_{k+1})\leq \eta (Q^\ast_k-Q_k(x_k))~~ \text{for given}~~ \eta\in(0,1].
\end{align}
Now, we list main assumptions involved in this paper as follows:
\begin{itemize}
	\item[A0]The set of minimizers of (\ref{model}) is nonempty, denoted by $\mathcal{X}$. The objective function $F$ attains its minimum $F^\ast=\inf_{x\in\mathcal{H}} F(x)\in\mathbb{R}$.
	\item[A1]$f,g:\mathcal{H}\rightarrow(-\infty,\infty]$ are two proper lower semi-continuous (lsc) convex functions with $\dom g\subset \interior(\dom f)$ and thus $\dom F=\dom g$.
	\item[A2] $f$ is Fr\'echet differentiable on an open set containing $\dom g$. Its gradient $\nabla f$ is continuous on $\dom g$.
	\item[A3]$\{\langle\cdot,\cdot\rangle_k\}$ is a sequence of inner products on $\mathcal{H}$, with induced norms $\{\|\cdot\|_k\}$ and associated positive operators $\{H_k\}$, i.e.,
	\begin{align}
	\forall k\geq0,~~~~H_k:\mathcal{H}\rightarrow\mathcal{H},~~~~\langle\cdot,\cdot\rangle_k=\langle\cdot,H_k\cdot\rangle.
	\end{align}
	There exist positive constants $M, m$ such that
	\begin{align}
	\forall k\geq0,~~~~~m\|\cdot\|^2\leq\|\cdot\|^2_k\leq M\|\cdot\|^2.
	\end{align}
	
	\item[A4]There exists an linearly convergent algorithm for the subproblem (\ref{sqa}) with a uniform parameter $\sigma$ for all $k\geq0$ such that
	\begin{align*}
	Q_k(y_l^{(k)})-Q_k^\ast\leq -(1-\sigma)^l Q_k^\ast, \forall l\geq 0
	\end{align*}
	where $0<\sigma<1,y_0^{(k)}:=x_k$.
\end{itemize}

The assumptions A0 and A1 are standard. The assumption A2 is weaker than the standard assumption that supposes $\nabla f$ to be Lipschitz continuous.
We do not assume any special structures on $H_k$ in A3. Note that the auxiliary function $Q_k$ is a regularized strongly convex function. The standard proximal gradient method, as shown in \cite[Theorem 2.1]{taylor2017exact}, could satisfy A4. This assumption is used to guarantee that each subproblem could be solved to satisfy the inexactness condition (\ref{inexact}) in a fixed number of iterations. Hence, the complexity of solving subproblems could not increase as $k\rightarrow +\infty$.

\section{The algorithm}
\subsection{Line Search}
For a directional line search method, we fix the direction $\bar x_{k+1}-x_k$, along which we search for a stepsize as large as possible. The main advantage of this kind of line search compared to \cite[Algorithm 2]{lee2018inexact} is that we only need to find an $\eta$-approximate minimizer of (\ref{sqa}) in each iteration.
And then we determine the next iterate by $x_{k+1}:=x_{k}+\alpha_k(\bar x_{k+1}-x_k)$. Below, we give several line search strategies to determine $\alpha_k$. For simplicity, we define
\begin{align*}
\Delta_k(x):=\langle \nabla f(x_k),x-x_k\rangle+g(x)-g(x_k).
\end{align*}

\noindent LS1. Let $\beta,\gamma,\in(0,1), \bar \alpha\in(0,1]$ and $\forall k\geq 0$,
\begin{align}
\alpha_k=\max\{&\alpha>0~|~\exists i\geq0,\alpha=\bar\alpha\beta^i,F(x_k+\alpha(\bar x_{k+1}-x_k))\nonumber\\
&-F(x_k)\leq \gamma \alpha\left(\langle\nabla f(x_k),\bar x_{k+1}-x_{k}\rangle\right.+g(\bar x_{k+1})-g(x_k))\}\label{ls1}.
\end{align}

\noindent LS2. Let $\beta\in(0,1),\gamma\in(0,m/2),\bar\alpha\in(0,1]$ and $\forall k\geq 0$,
\begin{align}
\alpha_k=\max&\{\alpha>0~|~\exists i\geq0,\alpha=\bar\alpha\beta^i,\alpha\|\nabla f(x_k+\alpha(\bar x_{k+1}-x_k))\nonumber\\
&-\nabla f(x_k)\|\leq \gamma\|x_{k+1}-x_{k}\|\}\label{ls2}.
\end{align}

\noindent LS3. Let $\beta,\gamma\in(0,1),\bar\alpha\in(0,1]$ and $\forall k\geq 0$,
\begin{align}
\alpha_k=&\max\left\{\alpha\right.>0~|~\exists i\geq0,\alpha=\bar\alpha\beta^i,F(x_k+\alpha(\bar x_{k+1}-x_k))-F(x_k)\nonumber\\
&\leq \alpha\left(\langle\nabla f(x_k),\bar x_{k+1}-x_{k}\rangle+g(\bar x_{k+1})-g(x_k)+\frac{\gamma}{2}\|\bar x_{k+1}-x_k\|_k^2\right)\}\label{ls3}.
\end{align}

\noindent LS4. Let $\beta,\gamma\in(0,1),\bar\alpha\in(0,1]$ and $\forall k\geq 0$,
\begin{align}
\alpha_k=&\max\left\{\alpha\right.>0~|~\exists i\geq0,\alpha=\bar\alpha\beta^i,f(x_{k+1})-f(x_k)\nonumber\\
&\leq \alpha\left(\langle\nabla f(x_k),\bar x_{k+1}-x_{k}\rangle+\frac{\gamma}{2}\|\bar x_{k+1}-x_k\|_k^2\right)\}\label{ls4}.
\end{align}

\subsection{Algorithmic Framework}
Now, we present the promised algorithmic framework of inexact SQA methods.

\begin{algorithm}[H]\label{framework}
	\SetAlgoLined

	\textbf{Initialization}: \\
	{\color{black}
		~~~~Given initial iterate $x_0\in\mathcal{H}$, $\eta\in(0,1]$\;}
	~~~~Choose $Linesearch(\cdot)$ from LS1-4 with proper parameters\;
	\For{$i=0,1,2,\cdots$}{
		Choose a symmetric $H_k$\;
		\textbf{Solving the subproblem inexactly}: \\
		~~~~Find an $\eta$-approximation $\bar x_{k+1}$ satisfying (\ref{inexact})\;
		$\alpha_k:=Linesearch(x_k,\bar x_{k+1})$\;
		$x_{k+1}:=x_k+\alpha_k(\bar x_{k+1}-x_k)$
	}
	\caption{Inexact Successive Quadratic Approximation with
		Linesearch}
\end{algorithm}

First of all, we state that the algorithmic framework is well defined. Its proof can be found in Appendix.
\begin{lemma}\label{welld}The stepsize $\alpha_k$ of LS1-4 exists.
\end{lemma}

The result above indicates that LS1-4 could find $\alpha_k$ by initializing $\alpha:=\bar \alpha$ and updating $\alpha:=\beta\alpha$ in finite algorithmic steps.

Next, we show a sufficient decrease property of the algorithmic framework.
Actually, we can derive that
\begin{align}
Q_k(\bar x_{k+1})=&g(\bar x_{k+1})-g(x_k)+\langle \nabla f(x_k),\bar x_{k+1}-x_k\rangle+\frac{1}{2}\|x_k-\bar x_{k+1}\|_k^2\label{Akd}\\
\leq &\eta\left(g(J_k)-g(x_k)+\langle \nabla f(x_k),J_k-x_k\rangle+\frac{1}{2}\|x_k-J_k\|_k^2\right)\nonumber\\
\leq &\eta\left(g(y)-g(x_k)+\langle \nabla f(x_k),y-x_k\rangle+\frac{1}{2}\|x_k-y\|_k^2-\frac{1}{2}\|y-J_k\|_k^2\right)\nonumber\\
= &\eta\left(g(x_k+\lambda (x-x_k))-g(x_k)+\langle \nabla f(x_k), \lambda (x-x_k)\rangle\right.\nonumber\\
&\left.+\frac{1}{2}\|\lambda (x-x_k)\|_k^2-\frac{1}{2}\|x_k+\lambda (x-x_k)-J_k\|_k^2\right)\nonumber\\
\leq &\eta \left(\lambda(F(x)-F(x_k))+\frac{\lambda^2}{2}\|x-x_k\|^2_k\right)\label{Ak}
\end{align}
for any $\lambda\in[0,1],x\in\mathcal{H}$ and $y:=x_k+\lambda (x-x_k)$, where the first inequality is due to that $\bar x_{k+1}$ is an $\eta$-approximate minimizer satisfying (\ref{inexact}), the second inequality follows from the strong convexity of $Q_k(\cdot)$, and the last inequality from the convexity of $F$. With this deduction, we have the following two results, whose proofs can be found in Appendix.
\begin{lemma}\label{FQ}With the proper parameters selected in LS1-4, we have
\begin{itemize}
  \item[(i)] $F(x_{k+1})-F(x_k)\leq \gamma\alpha_kQ_k(\bar x_{k+1})$ for LS1 and
  \item[(ii)] $F(x_{k+1})-F(x_k)\leq \alpha_kQ_k(\bar x_{k+1})$ for LS2-4.
\end{itemize}
\end{lemma}

\begin{lemma}\label{dec}
	For LS1-4, we have the sufficient decreasing property for all $k\geq 0$:
	\begin{align}
	F(x_{k+1})-F(x_k)\leq -\alpha_kc_1\|\bar x_{k+1}-x_k\|_k^2,
	\end{align}
	where $c_1$ is some positive constant. Thus,
	$\{F(x_k)\}$ is monotone decreasing and $\sum_{k=0}^{+\infty} \alpha_k\|\bar x_{k+1}-x_k\|^2< +\infty$.
\end{lemma}

\section{Complexity Analysis}
In this section, we will analyze the global complexity of the proposed algorithmic framework under mild assumptions. All proofs can be found
in Appendix.
\subsection{Linear Convergence Results}
In this subsection, we focus on convergence analysis under the quadratic growth condition. First, we introduce the optimal set strongly convexity condition (OSSC), which is presented in \cite{lee2018inexact} to get linear convergence. We say that a function $F$ satisfies OSSC if there exists $\mu>0$ such that for any $x\in\dom F$ and any $\lambda\in[0,1]$, it holds
\begin{align}\label{OSSC}
F(\lambda x+(1-\lambda)P_{\mathcal{X}}(x))\leq \lambda F(x)+(1-\lambda)F^\ast-\frac{\mu\lambda(1-\lambda)}{2}\|x-P_{\mathcal{X}}(x)\|^2,
\end{align}
where $P_{\mathcal{X}}(x):=\mathop{\arg\min}_{y\in\mathcal{X}}\|x-y\|$. Note that $\mathcal{X}$ is nonempty, convex and closed and hence $P_{\mathcal{X}}(x)$ is well-defined\cite[Theorem 1.2.3]{cegielski2012iterative}.

Below, we recall the quadratic growth condition.
\begin{definition}[\cite{zolezzi1977equiwellset,zhang2016new}] We say the function $F$ satisfies the $\mu$-quadratic growth (QG) condition if there exists $\mu>0$ such that
	\begin{align*}
	F(x)-F^\ast\geq\frac{\mu}{2}\|x-P_{\mathcal{X}}(x)\|^2,\forall x\in\dom F.
	\end{align*}
\end{definition}
Here, we claim that QG is strictly weaker than OSSC. For example, consider the function
\begin{align*}
F(x)=\left\{
\begin{array}{ll}
|x| &~~~~ \text{if}~~ |x|<1 \\
x^2 &~~~~ \text{else} \\
\end{array}
\right.,
\end{align*}
where $F(x)-F(0)\geq x^2$. It satisfies QG but not OSSC. Moreover, we observe that OSSC is sufficient for a nonsmooth extension of quasi strongly convexity, which is strictly stronger than QG; for details please refer to \cite{necoara2016linear}.

Now, we present the main result of this part.
\begin{theorem}\label{thm1}
	If $F$ satisfies the $\mu$-quadratic growth condition {\color{black}and $\inf_{k\geq0} \alpha_k\geq\underline{\alpha}$ for some $\underline{\alpha}>0$}. {\color{black}T}hen
	\begin{itemize}
		\item[(i)] The function value sequence $\{F(x_k)\}$ is Q-linearly convergent to $F^\ast$.
		\item[(ii)] The iterate sequence $\{x_k\}$ R-linearly converges to an optimal point $x^\ast$.
	\end{itemize}
\end{theorem}

Though the inexactness condition (\ref{inexact}) is hard to verify, we could use a fixed number of iterations $N_{inner}$. With the assumption A4, $N_{inner}$ iterations achieve (\ref{inexact}) with $\eta=1-(1-\sigma)^{N_{inner}}$.
Then we immediately have the following corollary.
\begin{corollary}\label{cor} Assume that the conditions in Theorem \ref{thm1} holds. Fixing inner iteration number $N_{inner}$ of the algorithm that satisfies A4, Algorithm 1 attains a solver $\tilde x$ such that
	\begin{align*}
	F(\tilde x)-F^\ast\leq \varepsilon
	\end{align*}	
	with $\mathcal{O}(\log(1/\varepsilon))$ inner iterations in total.
\end{corollary}

\subsection{Sublinear Convergence Results}
In this subsection, we drop the QG assumption of $F$. In order to illustrate the convergence of $\{F(x_k)\}$, we modify lemma in \cite{lee2018inexact} as follows
\begin{lemma}\label{sublemma}
	Assume we have three non-negative sequences $\{\delta_k\}_{k\geq 0}$, $\{\lambda_k\}_{k\geq 0}$, $\{A_k\}$ and a positive constant $c\in(0,1]$ such that
	\begin{align*}
	\delta_{k+1}\leq \delta_{k}+c\left(-\lambda_k\delta_k+\frac{A_k}{2}\lambda_k^2\right),~~~~\forall k\geq 0,\lambda_k\in[0,1].
	\end{align*}
	(i) If $A_k\leq \bar A$ for $k\geq 0$, where $\bar A$ is a positive constant, we have
	\begin{align}
	\delta_k\sim\mathcal{O}\left(k^{-1}\right).
	\end{align}
	(ii) If $\lim_{k\rightarrow 0} A_k=0$, we have
	\begin{align}
	\delta_k\sim o\left(k^{-1}\right).
	\end{align}
	
\end{lemma}

The following result improves the existing convergence rate of $\{F(x_k)\}_{k\geq 0}$ from $O(k^{-1})$ to $o(k^{-1})$. Denote
\begin{align*}
R_0:=\sup_{x:F(x)\leq F(x_0)}\|x-P_{\mathcal{X}}(x)\|,
\end{align*}
and assume $R_0$ to be finite. Since $F(x_k)\leq F(x_0)$ holds for $k\geq0$, we have
\begin{align}\label{r0}
\dist(x_k,\mathcal{X})\leq R_0,~~~~(k\geq0).
\end{align}

\begin{theorem}\label{generalConvergence}
	Suppose there exists $\underline\alpha>0$ such that $\alpha_k\geq\underline\alpha>0$ for all $k\in\mathbb{N}$ and $R_0$ is finite.  Then, $\{F(x_k)\}$ converges to $F^\ast$ sublinearly in the sense that
	\begin{align*}
	F(x_k)-F^\ast\sim\mathcal{O}\left(k^{-1}\right).
	\end{align*}
	Furthermore, if $\dist(x_k,\mathcal{X})\rightarrow 0$, then
	\begin{align*}
	F(x_k)-F^\ast\sim o\left(k^{-1}\right).
	\end{align*}
\end{theorem}

If the statement $[\forall\{x_k\}, F(x_k)\downarrow F^\ast\Rightarrow\dist(x_k,\mathcal{X})\rightarrow0]$ holds, the condition of the second conclusion in Theorem \ref{generalConvergence} will be automatically satisfied. Unfortunately, it is not true in general as the following counterexample illustrates:

Consider the function $F:\mathbb{R}^2\rightarrow (-\infty,\infty]$ satisfying lsc.,
\begin{align}
F(x,y)=\left\{
\begin{array}{cc}
x+\sqrt{x^2+y^2}&~~~~\text{if}~~x+y^2\leq1,\\
\infty&~~~~\text{otherwise}.
\end{array}
\right.
\end{align}
$F$ is convex since it is the sum of two convex functions. The optimal set is $\{(x,0)\in\mathbb{R}^2|x\leq0\}$ and the minimum is $0$. Consider the sequence $\{z_k\}\subset\dom F$ where for every $k\geq 0, z_k=(-\sum_{i=0}^k{1/(i+1)},1)$. It is obvious that $F(z_k)\downarrow 0$ but $\dist(z_k,\mathcal{X})\equiv1$.

Below, we propose several mild conditions, under which $F(x_k)\downarrow0$ implies $\dist(x_k,\mathcal{X})\rightarrow 0$.
\begin{proposition}\label{sufficient}
	If one of the following statements holds, then, any $\{x_k\}\subset\dom F$ satisfying $F(x_k)\downarrow F^\ast$ implies $\dist(x_k,\mathcal{X})\rightarrow 0$.
	\begin{enumerate}
		\item[(i)] The level set $C_0:=\{x\in\mathcal{H}|F(x)\leq F(x_0)\}$ is compact.
		\item[(ii)] $F$ is defined on $\mathbb{R}^n$ and its lineality space is equal to its recession cone(see definitions in \cite{bertsekas2009convex}), i.e.,
		\begin{align*}
		L_F=R_F.
		\end{align*}
		\item[(iii)]In particular, $F$ defined on $\mathbb{R}^n$ is level bounded, which implies
		\begin{align*} 
		L_F=R_F=\emptyset.
		\end{align*}
		
	\end{enumerate}
\end{proposition}

Note that for a globally $L$-smooth function $f$, there exists a positive number $\underline\alpha$ such that $\mathop{\lim\inf}_{k\rightarrow +\infty} \alpha_k= \underline\alpha> 0$. Therefore, using Theorem \ref{generalConvergence}, we have a slightly stronger convergence rate compared with \cite[Theorem 3]{lee2018inexact}.


Similar to Corollary \ref{cor}, we could use a fixed number of inner iterations in practice for general convex cases as well.

\begin{corollary}Suppose there exists $\underline\alpha>0$ such that $\alpha_k\geq\underline\alpha>0$ for all $k\in\mathbb{N}$ and A4 holds. With a fixed number of inner iterations replacing the stopping criterion (\ref{inexact}), Algorithm 5 attains a solver $\tilde x$ satisfying
	\begin{align*}
	F(\tilde x)-F^\ast\leq \varepsilon
	\end{align*}	
	with $\mathcal{O}(\varepsilon^{-1})$ inner iterations in total. Furthermore, the number of iterations is reduced to $o(\varepsilon^{-1})$ if $F$ is level bounded.
\end{corollary}

\subsection{Lower Bound for Stepsizes}
In this subsection, under a local gradient-Lispchitz-continuity condition, we prove that the stepsizes have a uniform positive lower bound, which guarantees that complexity of the line searches do not increase.
\begin{proposition}\label{dto0}If $\dist(x_k,\mathcal{X})\rightarrow 0$, $F(x_k)\downarrow F^\ast$ and $\nabla f$ is L-Lipschitz continuous on $\mathbb{B}_{\varepsilon}(\mathcal{X})\cap \dom F$ with $L>0,\varepsilon>0$ where
	\begin{align}
	\mathbb{B}_\varepsilon(\mathcal{X}):=\{x\in\mathcal{H}|\dist(x,\mathcal{X})\leq\varepsilon\}.
	\end{align}
	Denote $d_k:=\bar x_{k+1}-x_k$. Then
	\begin{align*}
	\lim_{k\rightarrow+\infty}\|d_k\|=0.
	\end{align*}
\end{proposition}

Equipped with the result above, the following lemma illustrates that stepsizes must have a uniform positive lower bound.

\begin{theorem}\label{thm2}Under the same conditions with Proposition \ref{dto0}, we have
	\begin{enumerate}
		\item[(i)] for LS1
		\begin{align}
		\liminf_{k\rightarrow+\infty}\alpha_k\geq \min\left\{1,\beta(1-\gamma)\frac{m(\eta+1+\sqrt{1-\eta})}{L(1+\sqrt{1-\eta})}\right\}.
		\end{align}
		\item[(ii)] For LS2
		\begin{align}\label{ls2inflim}
		\liminf_{k\rightarrow+\infty}\alpha_k\geq\min\left\{1,\frac{\beta\gamma}{L}\right\}.
		\end{align}
		\item[(iii)] For LS3 and LS4
		\begin{align}
		\liminf_{k\rightarrow+\infty}\alpha_k\geq\min\left\{1,\frac{\beta\gamma m}{L}\right\}.
		\end{align}
	\end{enumerate}
\end{theorem}

Now, from Theorems \ref{thm1}-\ref{thm2}, we can conclude that once the function value and iterate sequences are convergent, they must (sub)linearly converge under the local gradient-Lipschitz-continuity condition.

\section{Conclusion \& Future Work}
In this paper, we study the global complexity of an algorithmic framework of inexact SQA methods with four types of line search strategies under mild assumptions. On one hand, with the QG property and the uniform positive lower bound condition on stepsizes, we derive the Q-linear convergence of the function value sequence and the R-linear convergence of the iterate sequence. On the other hand,  without the QG property, we obtain the $o(k^{-1})$ complexity, which improves existing results. Finally, we give a uniform positive lower bound of the stepsizes for LS1-4 with the local gradient-Lipschitz-continuity assumption.

We believe that the new analysis developed in this paper might be extended to other related algorithms, such as inexact Bregman-type methods. We leave it as future work.



%
%

\bibliographystyle{plain}
\bibliography{ref}
\section*{Appendix}
\subsection*{A. Proof of Lemma \ref{welld}}
	It is easy to see that if $x_k=\bar x_{k+1}$, then the statement trivially holds. So we consider $x_k\neq \bar x_{k+1}$.
	If $x_k\in \mathcal{X}$, {\color{black}we have $Q_k^\ast=Q_k(x_k)$}, which implies that $Q_k(\bar x_{k+1})=Q_k(x_k)$. Due to the strong convexity of $Q_k(\cdot)$, it follows that $\bar x_{k+1}=x_k$. Therefore, we only need to consider $x_k\notin \mathcal{X}$, which implies {\color{black}$Q_k^\ast<0$} and hence $\Delta_k(\bar x_{k+1})<0$.

\noindent\textbf{LS1:}By contradiction suppose that for all $\alpha\in\mathcal{Q}:=\{\bar\alpha,\bar\alpha\beta,\bar\alpha\beta^2,\cdots\}$,
	\begin{align*}
	\frac{F(x_k+\alpha(\bar x_{k+1}-x_k))-F(x_k)}{\alpha}> \gamma\Delta_k(\bar x_{k+1}).
	\end{align*}
	With
	\begin{align}\label{gconvex}
	\alpha(g(\bar x_{k+1})-g(x_k)){\color{black}\geq} g(x_k+\alpha (\bar x_{k+1}-x_k))-g(x_k),
	\end{align} it follows that
	\begin{align*}
	\frac{f(x_k+\alpha(\bar x_{k+1}-x_k))-f(x_k)}{\alpha}+g(\bar x_{k+1})-g(x_k)> \gamma \Delta_k(\bar x_{k+1}).
	\end{align*}
	Taking $\alpha\downarrow 0$, due to $f$ is Fr\'echet differentiable at $x_k$, we obtain
	\begin{align*}
	\Delta_k(\bar x_{k+1})\geq\gamma\Delta_k(\bar x_{k+1}),
	\end{align*}
	a contradiction with $\Delta_k(\bar x_{k+1})<0$.
	
\noindent\textbf{LS2:} By contradiction suppose that for all $\alpha\in\mathcal{Q}:=\{\bar\alpha,\bar\alpha\beta,\bar\alpha\beta^2,\cdots\}$,
	\begin{align*}
	\|\nabla f(x_{k}+\alpha(\bar x_{k+1}-x_k))-\nabla f(x_k)\|>\gamma\|\bar x_{k+1}-x_k\|.
	\end{align*}
	Taking $\alpha\downarrow 0$, we have  $\nabla f(x_{k}+\alpha(\bar x_{k+1}-x_k))\rightarrow\nabla f(x_k)$. Then the continuity of $\nabla f$ at $x_k$ yields the contradiction $0\geq\gamma\|\bar x_{k+1}-x_k\|$.
	
\noindent\textbf{LS3:}  By contradiction suppose that for all $\alpha\in\mathcal{Q}:=\{\bar\alpha,\bar\alpha\beta,\bar\alpha\beta^2,\cdots\}$,
	\begin{align*}
	F(x_{k}+\alpha(\bar x_{k+1}-x_k))-F(x_k)>\alpha\left(\Delta_k(\bar x_{k+1})+\frac{\gamma}{2}\|\bar x_{k+1}-x_k\|^2_k\right).
	\end{align*}
	Using (\ref{gconvex}), dividing both sides by $\alpha$ and then taking $\alpha\downarrow 0$, due to that $f$ is Fr\'echet differentiable at $x_k$, we obtain
	\begin{align*}
	\Delta_k(\bar x_{k+1})\geq \Delta_k(\bar x_{k+1})+\frac{\gamma}{2}\|\bar x_{k+1}-x_k\|_k^2,
	\end{align*}
	a contradiction with $\bar x_{k+1}-x_k\neq 0$.

\noindent\textbf{LS4:}  By contradiction suppose that for all $\alpha\in\mathcal{Q}:=\{\bar\alpha,\bar\alpha\beta,\bar\alpha\beta^2,\cdots\}$,
	\begin{align*}
	f(x_{k}+\alpha(\bar x_{k+1}-x_k))-f(x_k)>\alpha\left(\langle\nabla f(x_k),\bar x_{k+1}-x_k\rangle+\frac{\gamma}{2}\|\bar x_{k+1}-x_k\|^2_k\right).
	\end{align*}
	Dividing both sides by $\alpha$ and then taking $\alpha\downarrow 0$, we obtain
	\begin{align*}
	\langle\nabla f(x_k),\bar x_{k+1}-x_k\rangle\geq \langle\nabla f(x_k),\bar x_{k+1}-x_k\rangle+\frac{\gamma}{2}\|\bar x_{k+1}-x_k\|^2,
	\end{align*}
	a contradiction with $\bar x_{k+1}-x_k\neq 0$.

\subsection*{B. Proof of Lemma \ref{FQ}}

	\textbf{LS1:} Combining (\ref{ls1}) with (\ref{Akd}), we obtain
	\begin{align*}
	\frac{1}{\alpha_k\gamma}\left(F(x_{k+1})-F(x_k)\right)+\frac{1}{2}\|\bar x_{k+1}-x_k\|^2_k\leq Q_k(\bar x_{k+1}),
	\end{align*}
	which implies the statement(i).\\
	\noindent\textbf{LS2:}
	Due to the convexity of $f$ and (\ref{ls2}), we have
	\begin{align}
	f(x_{k+1})-f(x_k)&\leq \langle\nabla f(x_{k+1}),x_{k+1}-x_k\rangle\nonumber\\
	&=\langle\nabla f(x_{k+1})-\nabla f(x_k),x_{k+1}-x_k\rangle+\langle \nabla f(x_k),x_{k+1}-x_k\rangle\nonumber\\
	&\leq \frac{\gamma}{\alpha_k}\|x_{k+1}-x_k\|^2+\langle \nabla f(x_k),x_{k+1}-x_k\rangle\nonumber\\
	&\leq \frac{\gamma}{\alpha_km}\|x_{k+1}-x_k\|_k^2+\langle \nabla f(x_k),x_{k+1}-x_k\rangle\nonumber\\
	&=\alpha_k\langle\nabla f(x_k),\bar x_{k+1}-x_k\rangle+\frac{\gamma\alpha_k}{m}\|\bar x_{k+1}-x_k\|_k^2.\label{p1}
	\end{align}
	Using the convexity of $g$ with $\alpha_k\in(0,1]$, we have
	\begin{align}\label{gk}
	g(x_{k+1})-g(x_k)\leq \alpha_k (g(\bar x_{k+1})-g(x_k)),
	\end{align}
	Adding (\ref{gk}) to (\ref{p1}) and then dividing $\alpha_k$ on both sides of the resulted inequality, we obtain
	\begin{align*}
	\frac{1}{\alpha_k}(F(x_{k+1})-F(x_k))+\left(\frac{1}{2}-\frac{\gamma}{m}\right)\|\bar x_{k+1}-x_k\|_k^2\leq Q_k(\bar x_{k+1}).
	\end{align*}
	Since we select $\gamma< m/2$, the statement(ii) for LS2 is proved.\\
	\noindent\textbf{LS3:} Combining (\ref{ls3}) with (\ref{Akd}), we obtain
	\begin{align}\label{qkineq1}
	\frac{1}{\alpha_k}\left(F(x_{k+1})-F(x_k)\right)+\frac{1-\gamma}{2}\|\bar x_{k+1}-x_{k}\|_k^2\leq Q_k(\bar x_{k+1}).
	\end{align}
	Since $\gamma\in(0,1)$, discarding the second term, then we proved the statement(ii) for LS3.\\
	\noindent\textbf{LS4:} Combining (\ref{ls4}) and (\ref{Akd}), we have
	\begin{align*}
	f(x_{k+1})-f(x_k)+\alpha_k(g(\bar x_{k+1})-g(x_k))\leq\alpha_kQ_k(\bar x_{k+1})+\alpha_k\left(\frac{\gamma-1}{2}\right)\|x_k-\bar x_{k+1}\|_k^2.
	\end{align*}
	Adding (\ref{gk}) to the inequality above and then dividing $\alpha_k$ on both sides of the resulted inequality, we obtain
	\begin{align*}
	\frac{1}{\alpha_k}\left(F(x_{k+1})-F(x_k)\right)+\frac{1-\gamma}{2}\|\bar x_{k+1}-x_k\|_k^2\leq Q_k(\bar x_{k+1}).
	\end{align*}
	Since $\gamma\in(0,1)$, the statement(ii) is proved.
\subsection*{C. Proof of Lemma \ref{dec}}
	Due to that $\bar x_{k+1}$ is an $\eta$-approximate minimizer, we have
	\begin{align}
	Q_k(\bar x_{k+1})&\leq\eta Q_k(J_k)\nonumber\\
	&\leq\eta Q_k(x_k+\lambda(\bar x_{k+1}-x_k))~~~~(\forall \lambda\in[0,1])\nonumber\\
	&\leq\eta\left(\lambda Q_k(\bar x_{k+1})-\frac{\lambda(1-\lambda)}{2}\|\bar x_{k+1}-x_k\|^2_k\right),\label{p2}
	\end{align}
	where the last inequality follows from the strong convexity of $Q_k$:
	\begin{align*}
	Q_k(x_k+\lambda(\bar x_{k+1}-x_k))\leq \lambda Q_k(x_k)+(1-\lambda)Q_k(x_k)-\frac{\lambda(1-\lambda)}{2}\|\bar x_{k+1}-x_k\|^2_k
	\end{align*}
	and the fact $Q_k(x_k)=0$.
	The inequality (\ref{p2}) leads to
	\begin{align*}
	Q_k(\bar x_{k+1})\leq-\frac{\eta\lambda(1-\lambda)}{2(1-\eta\lambda)}\|\bar x_{k+1}-x_k\|^2_k.
	\end{align*}
	Setting $\lambda:=\frac{1}{1+\sqrt{1-\eta}}$, which lies on $(1/2,1]$, we have:
	\begin{align}\label{Qkdk}
	Q_k(\bar x_{k+1})\leq-\frac{\eta}{2\left(1+\sqrt{1-\eta}\right)}\|\bar x_{k+1}-x_k\|_k^2.
	\end{align}
	Revoking Lemma \ref{FQ}, we obtain the following sufficient descent properties:
	\begin{align}
	\textbf{LS1}:&F(x_{k+1})-F(x_k)\leq -\frac{\alpha_k\gamma\eta}{2(1+\sqrt{1-\eta})}\|\bar x_{k+1}-x_k\|^2_k,\label{decforls1}\\
	\textbf{LS2-4}:&F(x_{k+1})-F(x_k)\leq -\frac{\alpha_k\eta}{2(1+\sqrt{1-\eta})}\|\bar x_{k+1}-x_k\|^2_k,\label{decforls2}
	\end{align}
	Therefore, $\{F(x_k)\}$ is monotone decreasing and
	\begin{align}\label{p3}
	\forall k\geq0, F(x_{k+1})-F(x_k)\leq-\alpha_k {\color{black}c_1}\|\bar x_{k+1}-x_k\|^2_k
	\end{align}
	for some positive constant ${\color{black}c_1}$. Summing up (\ref{p3}) for all $k\geq 0$, we have
	\begin{align*}
	\sum_{k=0}^{+\infty}\alpha_k\|\bar x_{k+1}-x_k\|_k^2\leq F(x_0)-\lim_{k\rightarrow \infty} F(x_k)\leq F(x_0)-F^\ast<\infty.
	\end{align*}

\subsection*{D. Proof of Lemma \ref{sublemma}}

	If $\delta_T=0$ for some $T$, then $\delta_k=0$ holds for $k\geq T$. Hence we assume $\delta_k\neq 0$ without loss of generality. The statement (i) is immediately obtained from \cite[Lemma 6]{lee2018inexact}. We consider the statement (ii). Since $A_k$ has a limit, then $A_k$ is upper bounded and thus (i) holds so that $\delta_{k}\rightarrow 0$.

	Since
	\begin{align*}
	\delta_{k+1}\leq\left\{
	\begin{array}{ll}
	(1-c)\delta_t+\frac{cA_k}{2} & \text{if~~} \delta_k> A_k\\
	\delta_k-\frac{c\delta_k^2}{2A_k} &\text{otherwise},
	\end{array}
	\right.
	\end{align*}
	and {\color{black}note that $\delta_k$ is monotone decreasing. With a slight abuse of notation that let $c/0=+\infty$ when $A_k=0$,} dividing $\delta_k\delta_{k+1}$ on both sides, then we obtain
	\begin{align*}
	\frac{1}{\delta_{k+1}}\geq\min\left\{\frac{1}{\delta_k}+\frac{1}{\delta_k}\frac{c}{2-c},\frac{1}{\delta_k}+\frac{c}{2A_k}\right\},
	\end{align*}
	which implies
	\begin{align*}
	\frac{1}{\delta_{k+1}}\geq\frac{1}{\delta_k}+z_k
	\end{align*}
	with $z_k=\min\{\frac{c}{\delta_k(2-c)},\frac{c}{2A_k}\}$. Since $z_k\rightarrow +\infty$, we immediately have
	\begin{align*}
	\delta_k\leq \frac{1}{\delta_0+\Sigma_{i=0}^{i=k-1}z_i}\sim o\left(\frac{1}{k}\right).
	\end{align*}

\subsection*{E. Proof of Proposition \ref{sufficient}}
	(i) Define a level set sequence $\{C_k\}$ associated with $\{l_k\}$, i.e., $\forall k\geq0$
	\begin{align*}
	C_k:=\{x\in\mathcal{H}|F(x)\leq l_k\},
	\end{align*}
	where we set $l_k:=F(x_k)(\forall k\geq0)$ and thus $l_k\downarrow F^\ast$. We will illustrate that $F(x_k)\rightarrow F^\ast$ implies $\dist(x,\mathcal{X})\rightarrow0$ by contradiction. Assume that there exist a subsequence $\{x_{k_i}\}$ and $D>0$ such that $\dist(x_{k_i},\mathcal{X})\geq D$ for every $i\geq 0$. The sequence $\{x_{k_i}\}$ is in the compact set $C_0$. Thus without any loss of generality, we assume $x_{k_i}\rightarrow x^\ast\in C_0$.
	Then we have $\dist(x^\ast,\mathcal{X})\geq D>0$.
	
	$C_k$ are closed due to $F$ is lsc. Because of the closedness of $C_k$ and $C_k\subset C_0$, each $C_k$ is compact, so is $\mathcal{X}$.

	$\forall \varepsilon>0, \exists I\geq 0$, $\forall i\geq I$, we have $\|x_{k_i}-x^\ast\|<\varepsilon/2$. Then we have $\mathbb{B}_{\varepsilon/2}(x^\ast)\cap C_{k_i}\neq \emptyset$ due to $x_{k_i}\in C_{k_i}$. Using $(\forall k\geq 0)C_{k+1}\subset C_k$, then we have $\forall k\geq 0, \mathbb{B}_{\varepsilon/2}(x^\ast)\cap C_{k}\neq \emptyset $. Denote $E_k:=\mathbb{B}_{\varepsilon/2}(x^\ast)\cap C_{k}$, then $E_k$ is compact due to compactness of $C_k$ and closedness of $\mathbb{B}_{\varepsilon/2}(x^\ast)$. Via $(\forall k\geq0)E_{k+1}\subset E_k$, we have $\cap_{k=0}^\infty E_k\neq\emptyset$, which leads to
	\begin{align*}
	\mathbb{B}_{\varepsilon/2}(x^\ast)\cap\mathcal{X}=\mathbb{B}_{\varepsilon/2}(x^\ast)\cap(\cap_{k=0}^\infty{C_k})=\cap_{k=0}^\infty E_k\neq\emptyset.
	\end{align*}
	Therefore, $x^\ast$ is in the closure of $\mathcal{X}$. Note that $\mathcal{X}$ is compact and hence closed. we have $x^\ast\in \mathcal{X}$, which contradicts $\dist(x^\ast,X)\geq D$.
	
	(ii) Let $C_k$ be defined as above. Then $C_k=L_F+\tilde{C}_k$, where $\tilde{C}_k\subset L_F^\perp$ and $\tilde {C}_k$ is compact\cite[Proposition 1.4.11]{bertsekas2009convex}. Each $x_k$ could be uniquely decomposed as $x_k=y_k+P_{L^\perp_F}(x_k)$ where $y_k\in L_F,P_{L^\perp_F}(x_k)\in\tilde{C}_k$. Define $\tilde F:\mathbb{R}^n\rightarrow (-\infty,\infty]$ as
	\begin{align*}
	\tilde{F}(x)=\left\{
	\begin{array}{cc}
	F(x)&~~\text{if~~}x\in L^\perp_F\\
	\infty &~~\text{otherwise}.
	\end{array}
	\right.
	\end{align*}
	Note that $\forall x\in\mathbb{R}^n$, we have $F(x)=F(P_{L^\perp_F}(x))=\tilde F(P_{L^\perp_F}(x))$. It is easy to show that $\tilde{F}$ is convex and lsc. Its minimum is $\inf_{x\in\mathbb{R}^n}\mathcal{\tilde F}(x)=F^\ast$ and the optimal set is $\tilde{ \mathcal{X}}=P_{L^\perp_F}(\mathcal{X})$. Consider the sequence $\{z_k\}$ where $z_k:=P_{L^\perp_F}(x_k)$. Then
	\begin{align}\label{limit}
	\tilde F(z_k)=F(x_k)\downarrow F^\ast = \inf_{x\in\mathbb{R}^n}\mathcal{\tilde F}(x).
	\end{align}
	Note that the set $\{x\in \mathbb{R}^n| \tilde F(x)\leq  \tilde F(z_0)\}=\{x\in L^\perp_F| F(x)\leq  F(x_0)\}=\tilde C_0$ is compact. Using (i) and (\ref{limit}), we have $\dist(z_k,\tilde{\mathcal{X}})\rightarrow 0$. Thus,
	\begin{align*}
	\dist(x_k,\mathcal{X})=\dist(x_k,\tilde{\mathcal{X}}+L_F)=\dist(z_k,\tilde{\mathcal{X}})\rightarrow 0
	\end{align*}
	which shows the statement (ii).

\subsection*{F. Proof of Proposition \ref{dto0}}

	Note that $Q_k(x_k)=0$ and $Q_k(\bar x_{k+1})\leq Q_k(x_k)$, then we have $Q_k(\bar x_{k+1})\leq 0$, i.e.,
	\begin{align*}
	\langle \nabla f(x_k),d_k\rangle+g(x_k+d_k)-g(x_k)+\frac{1}{2}\|d_k\|^2_k\leq0~~~(\forall k\geq 0){\color{black}.}
	\end{align*}
	Since $\|\cdot\|_k\geq m\|\cdot\|$, the inequality above immediately leads to
	\begin{align}\label{lessthan0}
	\langle \nabla f(x_k),d_k\rangle+g(x_k+d_k)-g(x_k)+\frac{m}{2}\|d_k\|^2\leq0~~~(\forall k\geq 0).
	\end{align}
	Suppose that there exists a subsequence $\{d_{k_i}\}$ and a positive number $D$ such that $\|d_{k_i}\|>D>0$ for all $i\geq 0$.
	For an arbitrary positive number $\delta$ which satisfies $0<\delta<\min\{\varepsilon/2,D\}$, since $\dist(x_k,\mathcal{X})\rightarrow 0$, we have $\dist(x_{k_i},\mathcal{X})<\delta$ for all large $i$. Denote  $\omega_{k_i}:=x_{k_i}+\delta d_{k_i}/\|d_{k_i}\|$. Then
	\begin{align*}
	\dist\left(\omega_{k_i},\mathcal{X}\right)&=\|\omega_{k_i}-P_\mathcal{X}(\omega_{k_i})\|=\left\|x_{k_i}+\delta\frac{ d_{k_i}}{\|d_{k_i}\|}-P_\mathcal{X}(x_{k_i})\right\|\\
	&\leq \delta+\|x_{k_i}-P_\mathcal{X}(x_{k_i})\|=\delta+\dist(x_{k_i},\mathcal{X})\\
	&<\delta+\dist(x_{k_i},\mathcal{X})\leq \varepsilon,
	\end{align*}
	which implies $\omega_{k_i}\in\mathbb{B}_\varepsilon(\mathcal{X})$ {\color{black}for all large $i$}. Also note that $\omega_{k_i}$ is on the line segment $[x_{k_i},x_{k_i}+d_{k_i}]\subset\dom F$, thus $\omega_{k_i}\in\mathbb{B}_{\varepsilon}(\mathcal{X})\cap \dom F$.
	The $L$-Lipschitz continuity of $\nabla f$ on $\mathbb{B}_\epsilon(\mathcal{X})\cap \dom F$ implies
	\begin{align}\label{locallips}
	\|\nabla f(x)-\nabla f(y)\|\leq L\|x-y\|~~~\forall x,y\in\mathbb{B}_\varepsilon(\mathcal{X})\cap \dom F.
	\end{align}
	Then using \cite[Lemma 2.64(i)]{bauschke2017convex},  we have
	\begin{align}\label{locallipineq}
	f(\omega_{k_i})-f(x_{k_i})\leq\frac{\delta}{\|d_{k_i}\|}\langle\nabla f(x_{k_i}),d_{k_i}\rangle+\frac{L\delta^2}{2}.
	\end{align}
	Due to the convexity of $g$ and $\delta/\|d_{k_i}\|< 1$, we have
	\begin{align*}
	g(\omega_{k_i})\leq \frac{\delta}{\|d_{k_i}\|}g(x_{k_i}+d_{k_i})+\left(1-\frac{\delta}{\|d_{k_i}\|}\right)g(x_{k_i}),
	\end{align*}
	a simple transformation of which yields
	\begin{align}\label{locallipg}
	g(\omega_{k_i})-g(x_{k_i})\leq\frac{\delta}{\|d_{k_i}\|}\left(g(x_{k_i}+d_{k_i})-g(x_{k_i})\right).
	\end{align}
	Combining (\ref{locallipineq}) and (\ref{locallipg}), {\color{black}for all large $i$,} we have
	\begin{align*}
	F\left(\omega_{k_i}\right)-F(x_{k_i})\leq \frac{\delta}{\|d_{k_i}\|}\left(\langle \nabla f(x_{\color{black}{k_i}}),d_{\color{black}{k_i}}\rangle+g(x_{\color{black}{k_i}}+d_{\color{black}{k_i}})-g(x_{\color{black}{k_i}})\right)+\frac{L\delta^2}{2}.
	\end{align*}
	Then employing (\ref{lessthan0}) leads to
	\begin{align}\label{preinflim}
	F\left(\omega_{k_i}\right)-F(x_{k_i})\leq& -\frac{\delta}{\|d_{k_i}\|}\frac{m\|d_{k_i}\|^2}{2}+\frac{L\delta^2}{2}
	\leq\frac{L}{2}\delta^2-\frac{\delta}{2}mD.
	\end{align}
	Since $F(x_k)\downarrow F^\ast$, taking the limit inferior on both sides of (\ref{preinflim}), we have
	\begin{align*}
	\liminf_{i\rightarrow +\infty} F(\omega_{k_i})-F^\ast\leq \frac{L}{2}\delta^2-\frac{\delta}{2}mD.
	\end{align*}
	By setting $\delta:=\min\{\varepsilon/4,mD/(2L),D/2\}$, we obtain
	\begin{align*}
	\liminf_{i\rightarrow +\infty}F(\omega_{k_i})< F^\ast,
	\end{align*}
	a contradiction with $F^\ast$ being the minimum. Therefore, we have $\|d_k\|\rightarrow 0$.

\subsection*{G. Proof of Theorem \ref{thm1}}

By setting $x:=P_{\mathcal{X}}(x_k)$ in (\ref{Ak}), we have
\begin{align}
Q_k(\bar x_{k+1})&\leq \eta \left(-\lambda (F(x_k)-F^\ast)+\frac{\lambda^2}{2}\|x_k-P_{\mathcal{X}}(x_k)\|_k^2\right)\label{qk}\\
&\leq \eta \left(-\lambda (F(x_k)-F^\ast)+\frac{M\lambda^2}{\mu}(F(x_k)-F^\ast)\right).\nonumber
\end{align}
The second inequality is due to the $\mu$-quadratic growth condition of $F$.
By setting $\lambda:=\min\{\mu/2M,1\}$, we have
\begin{align}\label{p4}
Q_k(\bar x_{k+1})\leq -\zeta\left(F(x_k)-F^\ast\right),
\end{align}
where $\zeta$ is a constant in $(0,1)$ satisfying
\begin{eqnarray*}
	\zeta=\left\{
	\begin{array}{ll}
		\eta\mu/(4M),&\text{if}~~\mu\leq 2M,\\
		1-M/\mu,&\text{else}.
	\end{array}
	\right.
\end{eqnarray*}
Note that $\alpha_k{\color{black}\geq \underline{\alpha}}$, using Lemma \ref{FQ}(i), together with (\ref{p4}) then we obtain
\begin{align*}
\forall k\geq 0,~~~~F(x_{k+1})-F^\ast\leq\left(1-\zeta{{\color{black}\underline \alpha}\gamma}\right)(F(x_k)-F^\ast)
\end{align*}
for LS1.
Similarly, for LS2-4 we have
\begin{align*}
\forall k\geq 0,~~~~F(x_{k+1})-F^\ast\leq\left(1-\zeta{{\color{black}\underline \alpha} }\right)(F(x_k)-F^\ast).
\end{align*}
Therefore, $\{F(x_k)\}$ is Q-linearly convergent to $F^\ast$:
\begin{align*}
F(x_k)-F^{\ast}\leq c^k (F(x_0)-F^\ast),
\end{align*}
where $c$ is a constant belonging to $(0,1)$.

We now prove that $\{x_k\}$ is R-linearly convergent. Using Lemma \ref{dec}, we have
	\begin{align*}
	\frac{c_1\|x_{k+1}-x_k\|^2}{\bar \alpha}&\leq \alpha_k c_1\|\bar x_{k+1}-x_k\|^2\\
	&\leq F(x_k)-F(x_{k+1})\leq F(x_k)-F^\ast\\
	&\leq c^k(F(x_0)-F^\ast).
	\end{align*}
	Thus, $\{\|x_{k+1}-x_k\|\}$ is linearly convergent satisfying
	\begin{align*}
	\forall k\geq0,~~~\|x_{k}-x_{k+1}\|\leq c^{\frac{k}{2}}\sqrt{\frac{\bar \alpha(F_0-F^\ast)}{c_1}},
	\end{align*}
	which implies that $\{x_k\}$ is a Cauchy sequence. By supposing $x_k\rightarrow x^\ast$, we have
	\begin{align*}
	\|x_k-x^\ast\|&\leq\sum_{i=k}^{+\infty}\|x_i-x_{i+1}\|\\
	&\leq \sqrt{\frac{\bar \alpha(F_0-F^\ast)}{c_1}}\sum_{i=k}^{+\infty}c^{\frac{i}{2}}\\
	&=\sqrt{\frac{\bar \alpha(F_0-F^\ast)}{c_1}}\frac{c^{\frac{k}{2}}}{1-\sqrt c}.
	\end{align*}
This is just the R-linear convergence of the iterate sequence $\{x_k\}$ and hence the proof is completed.

\subsection*{H. Proof of Theorem \ref{generalConvergence}}
Using (\ref{qk}), we have
\begin{align}\label{ineq1}
Q_k(\bar x_{k+1})&\leq \eta \left(-\lambda (F(x_k)-F^\ast)+\frac{\lambda^2M}{2}\dist(x_k,\mathcal{X})^2\right).
\end{align}
For LS1, combining (\ref{ineq1}) and Lemma \ref{FQ}(i), we have
\begin{align}\label{ineq2}
F(x_{k+1})-F(x_k)\leq\alpha_k\gamma\eta\left(\lambda(F^\ast-F(x_k))+\lambda^2\frac{M\dist(x_k,\mathcal{X})^2}{2}\right),\forall k\geq 0.
\end{align}
Let $\delta:=F(x_k)-F^\ast,c_k:=\alpha_k\gamma(1-\eta)$ and $\bar A:=MR_0^2$ in Lemma \ref{sublemma}; then we obtain
\begin{align*}
F(x_k)-F^\ast\leq\frac{MR_0^2+F(x_0)-F^\ast}{\sum_{i=0}^{k-1}a_i\gamma\eta}\sim\mathcal{O}\left(k^{-1}\right).
\end{align*}
Via the identical routine, similar results can also be obtained for LS2-4:
\begin{align*}
F(x_k)-F^\ast\leq\frac{MR_0^2+F(x_0)-F^\ast}{\sum_{i=0}^{k-1}a_i\eta}\sim \mathcal{O}\left(k^{-1}\right).
\end{align*}
Next, we will show $o(1/k)$ convergence in the function value sequence. Since $\dist(x_k,\mathcal{X})\rightarrow 0$, using (\ref{ineq2}) and Lemma \ref{sublemma}(ii), then for LS1-4, we have
\begin{align*}
F(x_k)-F^\ast\sim o\left(k^{-1}\right).
\end{align*}
The proof is completed.

\subsection*{I. Proof of Theorem \ref{thm2}}
	Since $\|d_k\|\rightarrow0$ and $\dist({\color{black}x_k},\mathcal{X})\rightarrow 0$, for all sufficiently large $k$, $\{x_k\}$ and $\{x_k+d_k\}$ will eventually fall into $\mathbb{B}_{\varepsilon}(\mathcal{X})\cap\dom F$. According to the $L$-Lipschitz continuity of $\nabla f$, we have
	\begin{align*}
	\|\nabla f(x_k+\alpha d_k)-\nabla f(x_k)\|\leq L\|\alpha d_k\|,~~~~\forall \alpha \in[0,1],
	\end{align*}
	which also implies
	\begin{align}\label{smallflips}
	f(x_k+\alpha d_k)-f(x_k)\leq \alpha \langle \nabla f(x_k),d_k\rangle+\frac{L\alpha^2}{2}\|d_k\|^2~~~~\forall \alpha\in[0,1].
	\end{align}
	
	(i)
	Adding (\ref{gk})(replacing $\alpha_k$ by $\alpha$) to (\ref{smallflips}),  for sufficiently large $k$ we have
	\begin{align}\label{flips}
	{\color{black}F(x_k+\alpha d_k)-F(x_k)\leq\alpha\Delta_k(\bar x_{k+1})+\frac{\alpha^2L}{2}\|d_k\|^2.}
	\end{align}
	From (\ref{Qkdk}), according to the relationship
	\begin{align} Q_k(\cdot)=\Delta_k(\cdot)+\frac{1}{2}\|\cdot\|_k^2\geq\Delta_k(\cdot)+\frac{m}{2}\|\cdot\|^2,
	\end{align}
	we have
	\begin{align}\label{deltad}
	\Delta_k(\bar x_{k+1})\leq-\left(\frac{m\eta}{2(1+\sqrt{1-\eta})}+\frac{m}{2}\right)\|d_k\|^2.
	\end{align}
	Combining (\ref{flips}) and (\ref{deltad}), canceling the term $\|d_k\|^2$, we have
	\begin{align}\label{ls1infa}
	F(x_k+\alpha d_k)-F(x_k)\leq \left(\alpha-\frac{\alpha^2L(1+\sqrt{1-\eta})}{m(\eta+1+\sqrt{1-\eta})}\right)\Delta_k(\bar x_{k+1}).
	\end{align}
	When $\alpha\leq\min\{1,(1-\gamma)\frac{m(\eta+1+\sqrt{1-\eta})}{L(1+\sqrt{1-\eta})}\}$, the stopping criterion in (\ref{ls1}) must hold. Hence (i) is proved.
	
	(ii) If $\alpha\leq\min\{1, \gamma/L\}$, {\color{black}then for all sufficiently large $k$, }
	\begin{align*}
	\alpha\|\nabla f(x_k+\alpha d_k)-\nabla f(x_k)\|&\leq \alpha L \|\alpha d_k\|\\
	&\leq \gamma \|x_{k+1}-x_k\|.
	\end{align*}
	which satisfies the stopping criterion in (\ref{ls2}) and thus (\ref{ls2inflim}) holds.
	
	(iii) Using (\ref{flips}), if $\alpha\leq\min\{1,\frac{\gamma m}{L}\}$, the stopping criterion in (\ref{ls3}) holds. The case of LS4 can be proved similarly.
\end{document}